\documentclass[12pt]{article}

\usepackage{amssymb}
\usepackage{amscd}

\def\carre {{\vrule height5pt width5pt depth0pt}}
\def\qed{{\hfill\carre\vskip1em}}

\newcommand{\email}[1]{{\footnote{#1}}}

\newcommand{\ad}{\mathop{\mathrm{ad}}}
\newcommand{\Aut}{\mathop{\mathrm{Aut}}}

\newcommand{\rg}{\mathop{\mathrm{rg}}}

\newcommand{\T}{\mathop{\mathrm{T}}}

\newcommand{\Sum}{\displaystyle\sum}

\begin{document}

\begin{center}{\huge
Indice et formes lin\'eaires stables dans 

les alg\`ebres de Lie
\vskip1em

\Large
Patrice Tauvel\email{E-mail: tauvel@mathlabo.univ-poitiers.fr} 
et
Rupert W.T. Yu\email{E-mail: yuyu@mathlabo.univ-poitiers.fr}
}
\vskip1em

UMR 6086 du C.N.R.S., 

D\'epartement de Math\'ematiques, 

Universit\'e de Poitiers, 

T\'el\'eport 2 -- BP 30179, 

Boulevard Marie et Pierre Curie, 

86962 Futuroscope Chasseneuil Cedex, 

FRANCE.
\end{center}

\noindent{\bf R\'esum\'e}
\vskip1em

\noindent On caract\'{e}rise, de mani\`{e}re
purement alg\'{e}brique, certaines formes lin\'{e}aires, dites
stables, introduites dans [7], sur une alg\`{e}bre de Lie. Comme
application, on d\'{e}termine l'indice d'une sous-alg\`{e}bre de Borel
d'une alg\`{e}bre de Lie semi-simple. On donne enfin un exemple d'une
sous-alg\`{e}bre parabolique d'une alg\`{e}bre de Lie semi-simple
n'admettant aucune forme lin\'{e}aire stable.
\vskip1em

\noindent{\bf Abstract}
\vskip1em

\noindent We characterize, in a purely algebraic manner, certain linear forms,
called stable, introduced in [7], on a Lie algebra. As an application,
we determine the index of a Borel subalgebra of a semi-simple Lie algebra.
Finally, we give an example of a parabolic subalgebra of a semi-simple Lie
algebra which does not admit any stable linear form.

\section{Formes lin\'{e}aires stables}

{\bf 1.1.}~Dans la suite, $\Bbbk$ est un corps commutatif
alg\'{e}briquement clos de caract\'{e}ristique nulle. Les alg\`{e}bres
de Lie consid\'{e}r\'{e}es sont d\'{e}finies et de dimension finie sur
$\Bbbk$. Si $V$ est un $\Bbbk$-espace vectoriel de dimension finie, on
le munit de la topologie de Zariski. Si $G$ est un groupe
alg\'{e}brique, on note $\mathcal{L}(G)$ son alg\`{e}bre de Lie. On
renvoie \`{a} [3] et [8] pour les concepts g\'{e}n\'{e}raux
utilis\'{e}s.

\smallskip

{\bf 1.2.}~Soient $\mathfrak{g}$ une alg\`{e}bre de Lie,
$\mathfrak{g}^{*}$ son dual. Si $\mathfrak{a}$ est un sous-espace
vectoriel de $\mathfrak{g}$, on note $\mathfrak{a}^{\bot}$ son
orthogonal dans $\mathfrak{g}^{*}$.

On fait op\'{e}rer $\mathfrak{g}$ dans $\mathfrak{g}^{*}$ au moyen de
la repr\'{e}sentation coadjointe. Ainsi, si $x,y \in \mathfrak{g}$ et
$f \in \mathfrak{g}^{*}$, on a :
$$
(x.f)(y) = f([y,x]). 
$$

Avec les notations pr\'{e}c\'{e}dentes, on d\'{e}finit une forme
bilin\'{e}aire altern\'{e}e $\Phi_{f}$ sur $\mathfrak{g}$ en posant :
$$
\Phi_{f}(x,y) = f([x,y]).
$$
Le noyau de $\Phi_{f}$ est not\'{e} $\mathfrak{g}^{f}$. On a :
$$
\big( \mathfrak{g}^{f}\big)^{\bot} = \mathfrak{g}.f \ , \ [\mathfrak{g},
\mathfrak{g}^{f}]^{\bot} = \{ g \in \mathfrak{g}^{*} ; \,
\mathfrak{g}^{f} \subset \mathfrak{g}^{g} \} .
$$

L'entier $\chi (\mathfrak{g}) = \inf \{ \dim \mathfrak{g}^{f}; \, f
\in \mathfrak{g}^{*} \}$ est appel\'{e} l'indice de $\mathfrak{g}$. On
dit que $f$ est r\'{e}gulier si $\dim \mathfrak{g}^{f} = \chi
(\mathfrak{g})$. L'ensemble $\mathfrak{g}_{r}^{*}$ des
\'{e}l\'{e}ments r\'{e}guliers de $\mathfrak{g}^{*}$ est un ouvert
non vide de $\mathfrak{g}^{*}$.

\smallskip

{\bf 1.3.}~Soient $\Aut \mathfrak{g}$ le groupe des automorphismes de
$\mathfrak{g}$ et $K$ un sous-groupe alg\'{e}brique de $\Aut
\mathfrak{g}$, d'alg\`{e}bre de Lie $\mathfrak{k}$. Le groupe $K$
op\`{e}re rationellement sur $\mathfrak{g}$ et
$\mathfrak{g}^{*}$. Pour $\alpha \in K$, $f \in \mathfrak{g}^{*}$ et
$x \in \mathfrak{g}$, on a $(\alpha.f)(x) = f\big(
\alpha^{-1}(x)\big)$. La diff\'{e}rentielle de cette op\'{e}ration
d\'{e}finit une action de $\mathfrak{k}$ sur $\mathfrak{g}^{*}$ ; si
$X \in \mathfrak{k}$, on a $(X.f)(x) = - f\big( X(x)\big)$. On pose :
$$
\mathfrak{k}_{f} = \{ X \in \mathfrak{k}; \, X.f = 0\} \ , \
\mathfrak{k}^{(f)} = \{ x \in \mathfrak{g}; \, (\mathfrak{k}.f)(x) = \{
0\} \} .
$$
On a ainsi $(\mathfrak{k}^{(f)})^{\bot} = \mathfrak{k}.f$. Si $K$
contient le groupe adjoint alg\'{e}brique de $\mathfrak{g}$, alors
$\ad_{\mathfrak{g}} \mathfrak{g} \subset \mathfrak{k}$, donc
$\mathfrak{k}^{(f)} \subset \mathfrak{g}^{f}$.

La forme bilin\'{e}aire $(X,x) \to f\big( X(x)\big)$ sur $\mathfrak{k}
\times \mathfrak{g}$ induit une forme bilin\'{e}aire non
d\'{e}g\'{e}n\'{e}r\'{e}e sur $(\mathfrak{k} / \mathfrak{k}_{f})
\times (\mathfrak{g} / \mathfrak{k}^{(f)})$. Par suite :
$$
\dim  \mathfrak{k} - \dim \mathfrak{k}_{f} = \dim
\mathfrak{g} - \dim \mathfrak{k}^{(f)} .
$$

{\bf D\'{e}finition.}~{\em On dit que $f \in \mathfrak{g}^{*}$ est
$K$-stable s'il existe un voisinage $V$ de $f$ dans $\mathfrak{g}^{*}$
tel que, pour tout $g \in V$, $\mathfrak{g}^{g}$ et $\mathfrak{g}^{f}$
soient $K$-conjugu\'{e}s. Si $K$ est le groupe adjoint alg\'{e}brique
de $\mathfrak{g}$ (resp. si $K = \Aut \mathfrak{g}$), une forme
lin\'{e}aire $K$-stable est dite stable (resp. faiblement stable).}

\smallskip

{\bf 1.4.~Lemme.}~{\em Soit $f \in \mathfrak{g}^{*}$. On suppose que
$K$ est connexe et que le morphisme
$$
\psi \colon K \times [\mathfrak{g}, \mathfrak{g}^{f}]^{\bot} \to
\mathfrak{g}^{*} \ , \ (\alpha , h) \to \alpha (h)
$$
est dominant. Alors :

{\em (i)} $f \in \mathfrak{g}_{r}^{*}$.

{\em (ii)} Soit $W = [\mathfrak{g}, \mathfrak{g}^{f}]^{\bot} \cap
\mathfrak{g}_{r}^{*}$. La restriction de $\psi$ \`{a} $K \times W$ est
une application ouverte.}

\smallskip

{\em Preuve.}~D'apr\`{e}s ([2], corollaire, p. 109), l'image de
$\psi$ contient un ouvert non vide de $\mathfrak{g}^{*}$ donc un
\'{e}l\'{e}ment $g \in \mathfrak{g}_{r}^{*}$. Si $\alpha \in K$ et $h
\in [\mathfrak{g}, \mathfrak{g}^{f}]^{\bot}$ v\'{e}rifient $g = \alpha
(h)$, il vient $\alpha \big( \mathfrak{g}^{f}\big) \subset \alpha
\big( \mathfrak{g}^{h}\big) = \mathfrak{g}^{g}$. D'o\`{u} $f \in
\mathfrak{g}_{r}^{*}$. 

(ii) Soit $\varphi \colon K \times W \to \mathfrak{g}^{*}$ le
morphisme induit par $\psi$. D'apr\`{e}s (i), $W$ est un ouvert non
vide de $[\mathfrak{g}, \mathfrak{g}^{f}]^{\bot}$, donc $\varphi$ est
dominant.

Soient $N$ le normalisateur de $\mathfrak{g}^{f}$ dans $K$, $\alpha
\in K$, et $g \in W$. On va d\'{e}terminer la fibre $S = \varphi^{-1}
\big( \varphi (\alpha , g)\big)$.

Soient $\beta \in K$ et $h \in W$ v\'{e}rifiant $\beta (h) = \alpha
(g)$, donc $h = \beta^{-1}\alpha (g)$. Comme $\mathfrak{g}^{h} =
\mathfrak{g}^{g} = \mathfrak{g}^{f}$, il vient $\beta^{-1}\alpha \in
N$, soit $\beta \in \alpha N$. Il est alors imm\'{e}diat que $S = \{
(\alpha \gamma, \gamma^{-1}(g)); \, \gamma \in N\}$.

Or, $S$ est le graphe du morphisme $\alpha N \to \mathfrak{g}^{*}$,
$\gamma \to \gamma^{-1}\alpha (g)$, donc est isomorphe \`{a} $N$. Comme
$N$ est une vari\'{e}t\'{e} pure, il r\'{e}sulte de ([2], corollaire
1, p. 113, et proposition 7, p. 141) que $\varphi$ est une application
ouverte. \qed

\smallskip

{\bf 1.5. Lemme.}~{\em Soient $K^{\circ}$ la composante neutre de $K$
et $f \in \mathfrak{g}^{*}$. Les conditions suivantes sont
\'{e}quivalentes :

{\em (i)} $f$ est $K^{\circ}$-stable.

{\em (ii)} $f$ est $K$-stable.

Si ces conditions sont v\'{e}rifi\'{e}es, on a $f \in
\mathfrak{g}_{r}^{*}$.} 

\smallskip

{\em Preuve.}~(i) $\Rightarrow$ (ii) C'est clair.

(ii) $\Rightarrow$ (i) Soit $U$ un voisinage ouvert de $f$ dans
$\mathfrak{g}^{*}$ tels que $\mathfrak{g}^{h}$ et $\mathfrak{g}^{f}$
soient $K$-conjugu\'{e}s pour tout $h \in U$. On a $U \cap
\mathfrak{g}_{r}^{*} \ne \emptyset$, donc $f \in
\mathfrak{g}_{r}^{*}$. Il en r\'{e}sulte que $W = [\mathfrak{g},
\mathfrak{g}^{f}]^{\bot} \cap \mathfrak{g}_{r}^{*}$ est un ouvert non
vide $[\mathfrak{g}, \mathfrak{g}^{f}]^{\bot}$ ; c'est l'ensemble des
$h \in \mathfrak{g}^{*}$ qui v\'{e}rifient $\mathfrak{g}^{h} =
\mathfrak{g}^{f}$. 

Si $g \in U$, il existe $\alpha \in K$ tel que $\mathfrak{g}^{g} =
\alpha \big( \mathfrak{g}^{f}\big)$, donc
$\mathfrak{g}^{\alpha^{-1}(g)} = \mathfrak{g}^{f}$, puis $g \in \alpha
(W)$. D\'{e}finissons le morphisme $\psi \colon K \times W \to
\mathfrak{g}^{*}$, $(\alpha , h) \to \alpha (h)$. Ce qui
pr\'{e}c\`{e}de montre que l'image $T$ de $\psi$ contient $U$.

Notons $K_{0} = K^{\circ}, K_{1} = \alpha_{1}K^{\circ}, \dots , K_{r} =
\alpha_{r}K^{\circ}$ les composantes irr\'{e}ductibles deux \`{a} deux
distinctes de $K$, avec $\alpha_{1}, \dots , \alpha_{r} \in K$. Posons
$S_{i} = K_{i} \times W$ et $T_{i} = \psi (S_{i})$. De $U \subset T_{0} \cup
\cdots \cup T_{r}$, on d\'{e}duit $\mathfrak{g}^{*} = \overline{T_{0}}
\cup \cdots \cup \overline{T_{r}}$. L'espace $\mathfrak{g}^{*}$
\'{e}tant irr\'{e}ductible, il existe un indice $i$ tel que
$\mathfrak{g}^{*} = \overline{T_{i}} = \overline{\alpha_{i}(T_{0})} =
\alpha_{i}(\overline{T_{0}})$. Par suite $\overline{T_{0}} =
\mathfrak{g}^{*}$, et le morphisme $\varphi \colon K^{\circ} \times W
\to \mathfrak{g}^{*}$ est dominant. D'apr\`{e}s 1.4, $\varphi$ est une
application ouverte. Alors $T_{0}$ est un voisinage ouvert de $f$ dans
$\mathfrak{g}^{*}$ et, pour $g \in T_{0}$, $\mathfrak{g}^{g}$ et
$\mathfrak{g}^{f}$ sont $K^{\circ}$-conjugu\'{e}s. \qed

\smallskip

{\bf 1.6.}~Si $\varphi \colon X \to Y$ est un morphisme de
vari\'{e}t\'{e}s, et si $x \in X$, on note $d \varphi_{x}$ la
diff\'{e}rentielle de $\varphi$ au point $x$ et $\T_{x}(X)$ l'espace
tangent de $X$ en $x$.

\smallskip

{\bf Lemme.}~{\em Soient $H$ un groupe alg\'{e}brique, $\mathfrak{h}$
son alg\`{e}bre de Lie, et $\pi$ une re\-pr\'{e}\-sen\-ta\-tion
rationnelle de $H$ dans un espace vectoriel $V$ de dimension finie. On
d\'{e}finit un morphisme de vari\'{e}t\'{e}s :
$$
\theta \colon H \times V \to V \ , \ (\beta , w) \to \pi (\beta )(w).
$$

{\em (i)} Si $(\alpha , v) \in H \times V$ et $(x, w) \in
\T_{\alpha}(H) \times V$, on a :
$$
d\theta_{(\alpha , v)}(x, w) = d \pi_{\alpha}(x)(v) + \pi (\alpha
)(w). 
$$

{\em (ii)} Soit $e$ l'\'{e}l\'{e}ment neutre de $H$. Pour tout
$(\alpha ,v) \in H \times V$, les applications lin\'{e}aires
$d\theta_{(\alpha , v)}$ et $d\theta_{(e,v)}$ ont m\^{e}me
rang.}

\smallskip

{\em Preuve.}~(i) D\'{e}finissons des morphismes
$$
p \colon H \to H \times V \ , \ \beta \to (\beta ,v) \ \mathrm{ et } \ q
\colon V \to H \times V \ , \ w \to (\alpha , w).
$$

De $\theta {\circ} p(\beta ) = \pi (\beta )(v)$ et $\theta {\circ} q
(w) = \pi (\alpha )(w)$, on d\'{e}duit :
$$
d (\theta {\circ}p)_{\alpha}(x) = d \pi_{\alpha}(x)(v) \ , \ d(\theta
{\circ}q)_{v}(w) = \pi (\alpha )(w).
$$
D'autre part, on a $dp_{\alpha}(x) = (x,0)$ et $dq_{v}(w) =
(0,w)$. Alors :
$$
\begin{array}{c}
d \pi_{\alpha}(x)(v) = d(\theta {\circ} p)_{\alpha} (x) =
d\theta_{(\alpha ,v)} {\circ} dp_{\alpha}(x) = d\theta_{(\alpha
,v)}(x,0), \\
\pi (\alpha )(w) = d(\theta {\circ} q)_{v}(w) = d \theta_{(\alpha ,v)}
{\circ} dq_{\alpha}(w) = d \theta_{(\alpha ,v)}(0, w).
\end{array}
$$
Comme $d \theta_{(\alpha ,v)}(x,w) = d\theta_{(\alpha , v)}(x, 0) + d
\theta_{(\alpha ,v)}(0, w)$, on a obtenu le r\'{e}sultat. 

(ii) Fixons $\alpha \in H$, et consid\'{e}rons les isomorphismes de
vari\'{e}t\'{e}s :
$$
\begin{array}{c}
\sigma \colon H \times V \to H \times V \ , \ (\beta , w) \to (\alpha
\beta , w) ,\\
\rho \colon V \to V \ , \ w \to \pi (\alpha )(w).
\end{array}
$$
On a $\theta {\circ} \sigma = \rho {\circ} \theta$. Si $v \in V$, il
vient donc $d(\theta {\circ} \sigma)_{(e,v)} = d (\rho {\circ} \theta
)_{(e,v)}$, c'est-\`{a}-dire $d \theta_{(\alpha ,v)} {\circ}
d\sigma_{(e,v)}= \pi (\alpha ) {\circ} d\theta_{(e,v)}$. D'o\`{u}
l'assertion. \qed

\smallskip

{\bf 1.7.~Th\'{e}or\`{e}me.}~{\em Soient $\mathfrak{g}$ une alg\`{e}bre
de Lie, $f \in \mathfrak{g}^{*}$, $K$ un sous-groupe alg\'{e}brique de
$\Aut \mathfrak{g}$, et $\mathfrak{k}$ l'alg\`{e}bre de Lie de
$K$. Les conditions suivantes sont \'{e}quivalentes :

{\em (i)} $[\mathfrak{g}, \mathfrak{g}^{f}] \cap \mathfrak{k}^{(f)} =
\{ 0\}$.

{\em (ii)} La forme lin\'{e}aire $f$ est $K$-stable.}

\smallskip

{\em Preuve.}~On peut supposer que $K$ est connexe (c'est clair pour
(i) $\Rightarrow$ (ii) et, pour (ii) $\Rightarrow$ (i), cela
r\'{e}sulte de 1.5). On note $e$ l'\'{e}l\'{e}ment neutre de $K$.

(i) $\Rightarrow$ (ii) On note $\psi \colon K \times [\mathfrak{g},
\mathfrak{g}^{f}]^{\bot} \to \mathfrak{g}^{*}$ le morphisme $(\alpha
, h) \to \alpha .h$. Soient $X \in \mathfrak{k}$ et $h \in
[\mathfrak{g}, \mathfrak{g}^{f}]^{\bot}$. D'apr\`{e}s 1.6, il
vient : 
$$
d \psi_{(e,f)}(X,h) = h + X.f.
$$
L'image de $d\psi_{(e,f)}$ est donc $[\mathfrak{g},
\mathfrak{g}^{f}]^{\bot} + \mathfrak{k}.f$, c'est-\`{a}-dire
$\mathfrak{g}^{*}$ d'apr\`{e}s l'hypoth\`{e}se. On sait alors que
$\psi$ est dominant. Soit $W = [\mathfrak{g}, \mathfrak{g}^{f}]^{\bot}
\cap \mathfrak{g}_{r}^{*}$. D'apr\`{e}s 1.4, $\psi (K \times W)$ est
un ouvert de $\mathfrak{g}^{*}$ contenant $f$ et, si $g \in \psi (K
\times W)$, $\mathfrak{g}^{f}$ et $\mathfrak{g}^{g}$ sont
$K$-conjugu\'{e}s.

(ii) $\Rightarrow$ (i) On a $f \in \mathfrak{g}_{r}^{*}$
(1.5), donc $W = [\mathfrak{g}, \mathfrak{g}^{f}]^{\bot} \cap
\mathfrak{g}_{r}^{*}$ est un ouvert non vide de
$[\mathfrak{g}, \mathfrak{g}^{f}]^{\bot}$ ; c'est l'ensemble des $h
\in \mathfrak{g}^{*}$ tels que $\mathfrak{g}^{h} =
\mathfrak{g}^{f}$.

Soit $\varphi \colon K \times W \to \mathfrak{g}^{*}$ le morphisme
$(\alpha ,h) \to \alpha .h$. Comme on l'a d\'{e}j\`{a} dit, $\varphi
(K \times W)$ est l'ensemble des $h \in \mathfrak{g}^{*}$ tels que
$\mathfrak{g}^{h}$ et $\mathfrak{g}^{f}$ soient
$K$-conjugu\'{e}s. 

Pour $g \in W$ et $(X,h) \in \mathfrak{k} \times
[\mathfrak{g}, \mathfrak{g}^{f}]^{\bot}$, il vient :
$$
d \varphi_{(e ,g)}(X,h) = h + X.g.
$$
L'image $T_{g}$ de $d \varphi_{(e ,g)}$ est donc $[\mathfrak{g},
\mathfrak{g}^{f}]^{\bot} + \mathfrak{k}.g$.

Notons $\mathfrak{n}$ le normalisateur de $\mathfrak{g}^{f}$ dans
$\mathfrak{k}$. On a :
$$
(X.g) ([\mathfrak{g}, \mathfrak{g}^{f}]) = g\big( [X(\mathfrak{g}),
\mathfrak{g}^{f}]\big) + g\big( [\mathfrak{g},
X(\mathfrak{g}^{f})]\big) = g\big( [\mathfrak{g},
X(\mathfrak{g}^{f})]\big) .
$$
Comme $\mathfrak{g}^{g} = \mathfrak{g}^{f}$, on en d\'{e}duit que
$X.g \in [\mathfrak{g}, \mathfrak{g}^{f}]^{\bot}$ si et seulement si
$X \in \mathfrak{n}$. D'o\`{u} :
$$
\dim (\mathfrak{k}.g \cap [\mathfrak{g}, \mathfrak{g}^{f}]^{\bot}) =
\dim \mathfrak{n}.g.
$$

D'autre part, si $X \in \mathfrak{k}_{g}$, on a $X \in \mathfrak{n}$,
car $X.g = 0$. D'o\`{u} :
$$
\dim \mathfrak{n}.g = \dim \mathfrak{n} - \dim \mathfrak{k}_{g} = \dim
\mathfrak{n} - \dim \mathfrak{k} + \dim \mathfrak{k}.g.
$$

On d\'{e}duit de tout ceci que :
$$
\dim T_{g} = \dim \mathfrak{k} - \dim \mathfrak{n} + \dim
([\mathfrak{g}, \mathfrak{g}^{f}]^{\bot}).
$$
En particulier, $\dim T_{g}$ ne d\'{e}pend pas de $g$. Il r\'{e}sulte
alors de 1.6 que, si $T_{(\alpha ,g)}$ est l'image de
$d\varphi_{(\alpha ,g)}$, $\dim T_{(\alpha ,g)}$ ne d\'{e}pend ni de
$\alpha$, ni de $g$. Or, d'apr\`{e}s les hypoth\`{e}ses, $\varphi$ est
dominant, donc il existe $(\alpha , g) \in K \times W$ tel que
$T_{(\alpha ,g)} = \mathfrak{g}^{*}$. On a alors $T_{(e,f)} =
\mathfrak{g}^{*}$, soit $[\mathfrak{g}, \mathfrak{g}^{f}]^{\bot} +
\mathfrak{k}.f = \mathfrak{g}^{*}$. Par suite, $[\mathfrak{g},
\mathfrak{g}^{f}] \cap \mathfrak{k}^{(f)} = \{ 0\}$. \qed

\smallskip

{\bf 1.8.~Corollaire.}~{\em Soit $f \in \mathfrak{g}^{*}$.

{\em (i)} Si $\mathfrak{g}$ est $\ad$-alg\'{e}brique, $f$ est
stable si et seulement si $[\mathfrak{g}, \mathfrak{g}^{f}] \cap
\mathfrak{g}^{f} = \{ 0\}$.

{\em (ii)} Si $[\mathfrak{g}, \mathfrak{g}^{f}] \cap
\mathfrak{g}^{f} = \{ 0\}$, $f$ est $K$-stable pour tout sous-groupe
alg\'{e}brique de $\Aut \mathfrak{g}$ v\'{e}rifiant
$\ad_{\mathfrak{g}} \mathfrak{g} \subset \mathcal{L}(K)$, donc en
particulier si $K$ contient le groupe adjoint alg\'{e}brique de
$\mathfrak{g}$. }

\smallskip

{\bf Remarque.}~Le r\'{e}sultat de 1.7 est \'{e}tabli dans [7]. On
n'a cependant utilis\'{e} ici que des m\'{e}thodes purement
alg\'{e}briques.

\section{Sous-alg\`{e}bres de Borel}

{\bf 2.1.}~Dans toute cette section, $\mathfrak{g}$ est une
alg\`{e}bre de Lie semi-simple, $G$ son groupe adjoint, $L$ sa forme
de Killing, $\mathfrak{h}$ une sous-alg\`{e}bre de Cartan de
$\mathfrak{g}$, et $R$ le syst\`{e}me de racines du couple
$(\mathfrak{g}, \mathfrak{h})$. Si $\alpha \in R$,
$\mathfrak{g}^{\alpha}$ est le sous-espace radiciel de $\mathfrak{g}$
associ\'{e} \`{a} $\alpha$. On fixe une base $\Pi$ de $R$, et on note
$R_{+}$ (resp. $R_{-}$) l'ensemble des racines positives
(resp. n\'{e}gatives) correspondant. On pose :
$$
\mathfrak{n} = \Sum_{\alpha \in
R_{+}} \mathfrak{g}^{\alpha} \ , \ \mathfrak{b} = \mathfrak{h} \oplus
\mathfrak{n} .
$$

Soient $\alpha , \beta \in R$. On fixe $X_{\alpha} \in
\mathfrak{g}^{\alpha} \backslash \{ 0\}$, et on note $H_{\alpha}$
l'unique \'{e}l\'{e}ment de $[\mathfrak{g}^{\alpha},
\mathfrak{g}^{-\alpha}]$ tel que $\alpha (H_{\alpha}) = 2$. On a
$\beta (H_{\alpha}) \in \mathbb{Z}$. Si $\lambda \in
\mathfrak{h}^{*}$, on \'{e}crit $\langle \lambda , \alpha^{\vee}
\rangle$ pour $\lambda (H_{\alpha})$.

Pour toute partie $S$ de $\Pi$, on note $\mathbb{Z}S$
(resp. $\mathbb{N}S$) l'ensemble des combinaisons lin\'{e}aires a
coefficients entiers (resp. \`{a} coefficients entiers positifs ou
nuls) des \'{e}l\'{e}ments de $S$. On pose :
$$
R^{S} = R \cap \mathbb{Z}S \ , \ R_{+}^{S} = R \cap \mathbb{N}S =
R_{+} \cap R^{S}.
$$
L'ensemble $R^{S}$ est un syst\`{e}me de racines dans le sous-espace
de $\mathfrak{h}^{*}$ qu'il engendre, $S$ est une base de $R^{S}$, et
$R_{+}^{S}$ est l'ensemble des racines positives correspondant. Si $S$
est une partie connexe de $\Pi$, le syst\`{e}me de racines $R^{S}$ est
irr\'{e}ductible, et on note $\varepsilon_{S}$ la plus grande racine
de $R^{S}$.

Supposons $S$ connexe. Pour toute racine $\alpha \in R_{+}^{S}
\backslash \{ \varepsilon_{S} \}$, on a $\langle \alpha ,
\varepsilon_{S}^{\vee} \rangle \in \{ 0,1\}$. Si $T$ est l'ensemble
des \'{e}l\'{e}ments $\alpha$ de $R^{S}$ qui v\'{e}rifient $\langle
\alpha , \varepsilon_{S}^{\vee}\rangle = 0$, on voit que $T$ est un
syst\`{e}me de racines dans le sous-espace de $\mathfrak{h}^{*}$ qu'il
engendre, et que $\{ \alpha \in S ; \, \langle \alpha ,
\varepsilon_{S}^{\vee} \rangle = 0\}$ est une base de $T$.

Si $\alpha \in T \cap R_{+}^{S}$, on a $\alpha + \varepsilon_{S}
\notin R$ et $\alpha - \varepsilon_{S} \notin R$ ([1], proposition 9,
p. 149). Ainsi, si $\alpha \ne \varepsilon_{S}$, les racines $\alpha$
et $\varepsilon_{S}$ sont fortement orthogonales.

\smallskip

{\bf 2.2.}~On va rappeler la construction et quelques
propri\'{e}t\'{e}s d'un ensemble de racines deux \`{a} deux fortement
orthogonales dans $R$ (voir [5] et [6]).

Soit $S$ une partie de $\Pi$. Par r\'{e}currence sur le cardinal de
$S$, on d\'{e}finit un sous-ensemble $\mathcal{K}(S)$ de
l'ensemble des parties de $S$ de la mani\`{e}re suivante :

a) $\mathcal{K}(\emptyset ) = \emptyset$.

b) Si $S_{1}, \dots , S_{r}$ sont les composantes connexes de $S$,
on a :
$$
\mathcal{K}(S) = \mathcal{K}(S_{1}) \cup \cdots \cup
\mathcal{K}(S_{r}) .
$$

c) Si $S$ est connexe, alors :
$$
\mathcal{K}(S) = \{ S\} \cup \mathcal{K} (\{ \alpha \in S; \, \langle
\alpha , \varepsilon_{S}^{\vee}\rangle = 0\} ).
$$

On a le r\'{e}sultat suivant :

{\bf Lemme.}~{\em {\em (i)} Tout $K \in \mathcal{K}(S)$ est une partie
connexe de $\Pi$.

{\em (ii)} Si $K, K' \in \mathcal{K}(S)$, alors ou $K \subset K'$, ou
$K' \subset K$, ou $K$ et $K'$ sont des parties disjointes de $S$
telles que $\alpha + \beta \notin R$ pour $\alpha \in R^{K}$ et $\beta
\in R^{K'}$.

{\em (iii)} Si $K$ et $K'$ sont des \'{e}l\'{e}ments distincts de
$\mathcal{K}(S)$, $\varepsilon_{K}$ et $\varepsilon_{K'}$ sont
fortement orthogonales.}

\smallskip

{\bf 2.3.}~Conservons les notations pr\'{e}c\'{e}dentes. Dans la
suite, on note $\mathcal{K}$ pour $\mathcal{K}(\Pi )$ et, si $K \in
\mathcal{K}$, on pose :
$$
\Gamma^{K} = \{ \alpha \in R^{K}; \, \langle \alpha ,
\varepsilon_{K}^{\vee}\rangle > 0\} \ , \ \Gamma_{0}^{K} = \Gamma^{K}
\backslash \{ \varepsilon_{K}\} \ , \ \mathfrak{a}_{K} =
\Sum_{\alpha \in \Gamma^{K}}
\mathfrak{g}^{\alpha}.
$$

{\bf Lemme.}~{\em Soient $K,K' \in \mathcal{K}$, $\alpha , \beta \in
\Gamma^{K}$, et $\gamma \in \Gamma^{K'}$.

{\em (i)} On a $\Gamma^{K} = R_{+}^{K} \backslash \{ \delta \in
R_{+}^{K}; \, \langle \delta , \varepsilon_{K}^{\vee} \rangle = 0 \}$.

{\em (ii)} L'ensemble $R_{+}$ est la r\'{e}union disjointe des
$\Gamma^{K''}$ pour $K'' \in \mathcal{K}$, et
$\mathfrak{a}_{K''}$ est une alg\`{e}bre de Heisenberg, de
centre $\mathfrak{g}^{\varepsilon_{K''}}$.

{\em (iii)} Si $\alpha + \beta \in R$, on a $\alpha + \beta =
\varepsilon_{K}$.

{\em (iv)} Si $\alpha + \gamma \in R$, alors ou $K \subset K'$ et
$\alpha + \gamma \in \Gamma^{K'}$, ou $K' \subset K$ et $\alpha +
\gamma \in \Gamma^{K}$. }

\smallskip

{\bf 2.4.}~Il est clair que le cardinal de $\mathcal{K}(\Pi )$ ne
d\'{e}pend que de $\mathfrak{g}$, mais  pas de $\mathfrak{h}$ ou de
$\Pi$. On note cet entier $k_{\mathfrak{g}}$. On donne dans le tableau
suivant la valeur de $k_{\mathfrak{g}}$ pour les diff\'{e}rents types
d'alg\`{e}bres de Lie simples. Si $r$ est un nombre rationnel, $[r]$
d\'{e}signe sa partie enti\`{e}re.

\medskip

\begin{center}
\begin{tabular}{|c|c|c|c|c|c|c|c|c|c|}\hline
\strut\rule[-1ex]{0ex}{4ex} & $A_{\ell}, \ell \geqslant 1$ &
$B_{\ell}, \ell \geqslant 2$ & $C_{\ell}, \ell \geqslant 3$ &
$D_{\ell} , \ell \geqslant 4$ & $E_{6}$ & $E_{7}$ & $E_{8}$ & $F_{4}$
& $G_{2}$ \\ \hline
$k_{\mathfrak{g}}$ \strut\rule[-2.8ex]{0ex}{7ex}& $\left[
\displaystyle\frac{\ell+1}{2}\right]$ & $\ell$ & $\ell$ & $2 \left[
\displaystyle\frac{\ell}{2}\right]$& $4$ & $7$ & $8$ & $4$ & $2$ \\ \hline
\end{tabular}
\end{center}

\smallskip

{\bf 2.5.~Lemme.}~{\em Soient $x \in \mathfrak{b}$ et
$$
u = \Sum_{K \in \mathcal{K}}
X_{-\varepsilon_{K}}. 
$$
Les conditions suivantes sont \'{e}quivalentes :

{\em (i)} $x \in \mathfrak{h}$ et $\varepsilon_{K}(x) = 0$ pour tout
$K \in \mathcal{K}$.

{\em (ii)} $[x,u] \in \mathfrak{n}$.}

\smallskip

{\em Preuve.}~(i) $\Rightarrow$ (ii) Si (i) est v\'{e}rifi\'{e}, on a
$[x,u] = 0$.

(ii) $\Rightarrow$ (i) Ecrivons
$$
x = h + \Sum_{\alpha \in R_{+}}
a_{\alpha} X_{\alpha},
$$
avec $h \in \mathfrak{h}$ et $a_{\alpha} \in \Bbbk$ pour $\alpha \in
R_{+}$. On obtient 
$$
[x,u] \in \Sum_{K \in \mathcal{K}}
a_{\varepsilon_{K}} \lambda_{K} H_{\varepsilon_{K}} + \mathfrak{n} +
\mathfrak{n}_{-}, 
$$
o\`{u} $\lambda_{K} \in \Bbbk \backslash \{ 0\}$ pour tout $K \in
\mathcal{K}$. Les $\varepsilon_{K}$ \'{e}tant deux \`{a} deux
fortement orthogonales, les $H_{\varepsilon_{K}}$ sont
lin\'{e}airement ind\'{e}pendants. Par suite, $a_{\varepsilon_{K}} =
0$ pour tout $K \in \mathcal{K}$.

Soient $K \in \mathcal{K}$, $\alpha \in \Gamma_{0}^{K}$, et $\beta =
\varepsilon_{K} - \alpha$. Il vient
$$
\begin{array}{rl}
[x,u] & = \lambda a_{\alpha} X_{-\beta} + a_{\alpha}
\Sum_{K' \in \mathcal{K}
\backslash \{ K\} } [X_{\alpha} ,X_{-\varepsilon_{K'}}] \\
{} & \strut\phantom{abcd} +
\Sum_{K' \in \mathcal{K},\gamma
\in R_{+} \backslash \{ \alpha \}} a_{\gamma} [X_{\gamma},
X_{-\varepsilon_{K'}}] -\Sum_{K'
\in \mathcal{K}} \varepsilon_{K'}(h) X_{-\varepsilon_{K'}},
\end{array}
$$
avec $\lambda \in \Bbbk \backslash \{ 0\}$.

On a $\beta \ne \varepsilon_{K}$ et $\beta \ne \varepsilon_{K'}$ si $K
\ne K'$, car $\langle \beta ,\varepsilon_{K}^{\vee} \rangle = 1$ et
$\langle \varepsilon_{K'}, \varepsilon_{K}^{\vee} \rangle = 0$ si $K
\ne K'$. 

Supposons $a_{\alpha} \ne 0$. Il existe n\'{e}cessairement $K' \in
\mathcal{K}$ et $\gamma \in R_{+} \backslash \{ \alpha \}$ tels que
$\beta = \varepsilon_{K'} - \gamma$. 

On a $K' \ne K$, car $\mathfrak{a}_{K}$ est une alg\`{e}bre de
Heisenberg. Soit $K'' \in \mathcal{K}$ tel que $\gamma \in
\Gamma^{K''}$. A nouveau, $K'' \ne K$. D'autre part, comme $\beta +
\gamma = \varepsilon_{K'} \in R$, on d\'{e}duit que $K \subset K''$,
ou $K'' \subset K$ (2.3).

Supposons $K'' \subset K$. Les racines $\varepsilon_{K'}$ et $\gamma$
sont fortement orthogonales \`{a} $\varepsilon_{K}$. D'o\`{u} :
$$
1 = \langle \beta , \varepsilon_{K}^{\vee} \rangle = \langle
\varepsilon_{K'} - \gamma , \varepsilon_{K}^{\vee} \rangle = 0.
$$
Contradiction.

Supposons $K \subset K''$. On a $\gamma \ne \varepsilon_{K''}$, car
$\varepsilon_{K''}$ est la plus grande racine de $R^{K''}$. Ainsi :
$$
\langle \gamma , \varepsilon_{K''}^{\vee} \rangle = 1 \ , \ \langle
\beta , \varepsilon_{K''}^{\vee} \rangle = 0.
$$
Par suite :
$$
1 = \langle \gamma , \varepsilon_{K''}^{\vee} \rangle = \langle \beta
+ \gamma , \varepsilon_{K''}^{\vee} \rangle = \langle
\varepsilon_{K'}, \varepsilon_{K''}^{\vee}\rangle \in \{ 0,2\} .
$$
C'est absurde.

On a donc prouv\'{e} que $a_{\alpha} =0$ pour tout $\alpha \in
R_{+}$. Il en r\'{e}sulte que $x \in \mathfrak{h}$, et que :
$$
[x,u] = - \Sum_{K \in \mathcal{K}}
\varepsilon_{K}(h) X_{-\varepsilon_{K}}. 
$$
D'o\`{u} le r\'{e}sultat. \qed

\smallskip

{\bf 2.6.~Lemme.}~{\em Soient $\mathfrak{a}$ une sous-alg\`{e}bre de
Lie d'une alg\`{e}bre de Lie semi-simple $\mathfrak{g}$, et $f \in
\mathfrak{a}^{*}$. On suppose que $\mathfrak{a}^{f}$ est une
alg\`{e}bre de Lie commutative compos\'{e}e d'\'{e}l\'{e}ments
semi-simples. Alors $[\mathfrak{a}, \mathfrak{a}^{f}] \cap
\mathfrak{a}^{f} = \{ 0\}$.}

\smallskip

{\em Preuve.}~D'apr\`{e}s les hypoth\`{e}ses, il existe un sous-espace
vectoriel $\mathfrak{r}$ de $\mathfrak{a}$ v\'{e}rifiant $\mathfrak{a}
= \mathfrak{a}^{f} \oplus \mathfrak{r}$ et $[\mathfrak{a}^{f},
\mathfrak{r}] \subset \mathfrak{r}$.  On en d\'{e}duit $[\mathfrak{a},
\mathfrak{a}^{f}] \subset \mathfrak{r}$. D'o\`{u} l'assertion. \qed

\smallskip

{\bf 2.7.~Th\'{e}or\`{e}me.}~{\em Soient $\mathfrak{g}$ une
alg\`{e}bre de Lie semi-simple, $\rg (\mathfrak{g})$ son rang,
$k_{\mathfrak{g}}$ l'entier de {\em 2.4}, et $\mathfrak{b}$ une
sous-alg\`{e}bre de Borel de $\mathfrak{g}$.

{\em (i)} L'espace $\mathfrak{b}^{*}$ contient un \'{e}l\'{e}ment
stable.

{\em (ii)} L'indice de $\mathfrak{b}$ est \'{e}gal \`{a} $\rg
(\mathfrak{g}) - k_{\mathfrak{g}}$.}

\smallskip

{\em Preuve.}~Si $\mathfrak{a}$ est une sous-alg\`{e}bre de Lie de
$\mathfrak{g}$, et si $y \in \mathfrak{g}$, on d\'{e}finit une forme
lin\'{e}aire $\varphi_{\mathfrak{a}}(x)$ sur $\mathfrak{a}$ en posant,
pour $x \in \mathfrak{a}$ :
$$
\varphi_{\mathfrak{a}}(y) (x) = L(y,x).
$$

Soient $u$ comme en 2.5 et $f = \varphi_{\mathfrak{b}}(u) \in
\mathfrak{b}^{*}$. On a $\mathfrak{b}^{f} = \{ x \in \mathfrak{b}; \,
[x,u] \in \mathfrak{n}\}$. Compte tenu de 2.5 et 2.6, on voit que $f$
est un \'{e}l\'{e}ment stable de $\mathfrak{b}^{*}$. L'assertion (ii)
r\'{e}sulte alors de 2.5, car les $\varepsilon_{K}$, pour $K \in
\mathcal{K}$, sont deux \`{a} deux fortement orthogonales, donc
lin\'{e}airement ind\'{e}pendantes. \qed

\smallskip

{\bf Remarque.}~L'assertion de 2.7, (ii) est obtenue sous une forme
diff\'{e}rente en ([9], theorem 4.1).

\section{Un contre-exemple}

{\bf 3.1.}~Il est donn\'{e} dans [7] des exemples d'alg\`{e}bres de
Lie ne poss\'{e}dant aucune forme lin\'{e}aire faiblement stable.

Soient $\mathfrak{g}$ une alg\`{e}bre de Lie semi-simple et
$\mathfrak{p}$ une sous-alg\`{e}bre parabolique de $\mathfrak{g}$. Il
est affirm\'{e} dans [4] que $\mathfrak{p}^{*}$ contient un
\'{e}l\'{e}ment stable. Nous allons donner un exemple prouvant que
cette affirmation est inexacte.

\smallskip

{\bf 3.2.}~On conserve les notations et conventions de 2.1, et on
suppose que $\mathfrak{g}$ est simple de type $D_{4}$. On note $\Pi =
\{ \alpha_{1}, \alpha_{2}, \alpha_{3}, \alpha_{4}\}$, la racine
$\alpha_{2}$ \'{e}tant li\'{e}e aux $\alpha_{i}$ pour $i \ne 2$ (c'est
la num\'{e}rotation des racines de [1], p. 256). On note $G$ le
groupe adjoint de $\mathfrak{g}$, et on pose $\mathfrak{p} = \Bbbk
X_{-\alpha_{2}} + \mathfrak{b}$. Ainsi, $\mathfrak{p}$ est une
sous-alg\`{e}bre parabolique minimale de $\mathfrak{g}$. On sait que
$\ad_{\mathfrak{g}} \mathfrak{b}$ et $\ad_{\mathfrak{g}} \mathfrak{p}$
sont des sous-alg\`{e}bres de Lie alg\'{e}briques de
$\ad_{\mathfrak{g}} \mathfrak{g}$. Soit $B$ (resp. $P$) le plus petit
sous-groupe alg\'{e}brique de $G$ tel que $\mathcal{L}(B) =
\ad_{\mathfrak{g}} \mathfrak{b}$ (resp. $\mathcal{L}(P) =
\ad_{\mathfrak{g}} \mathfrak{p}$). On a $B \subset P$, et $B$
(resp. $P$) s'identifie au groupe adjoint de $\mathfrak{b}$
(resp. $\mathfrak{p}$). 

L'ensemble $\mathcal{K}$ est constitu\'{e} de $\Pi$, $\alpha_{1},
\alpha_{3}, \alpha_{4}$. On en d\'{e}duit que l'ensemble $\mathcal{E}$
des $\varepsilon_{K}$, $K \in \mathcal{K}$, est :
$$
\mathcal{E} = \{ X_{\alpha_{1}}, X_{\alpha_{3}}, X_{\alpha_{4}},
X_{\alpha_{1} + 2 \alpha_{2} + \alpha_{3} + \alpha_{4}} \} .
$$

On conserve la notation $\varphi_{\mathfrak{a}}$ de la preuve de
2.7. Posons :
$$
u = X_{-\alpha_{1}} + X_{-\alpha_{3}} + X_{-\alpha_{4}} +
X_{-\alpha_{1} - 2 \alpha_{2} - \alpha_{3} - \alpha_{4}} \ , \ g =
\varphi_{\mathfrak{b}}(u). 
$$
D'apr\`{e}s 2.7, $g$ est un \'{e}l\'{e}ment stable de
$\mathfrak{b}^{*}$. D'autre part, l'indice de $\mathfrak{b}$ est nul
(2.4 et 2.7). Il en r\'{e}sulte que $\mathfrak{b}^{*}$ contient une
$B$-orbite ouverte. C'est celle de $g$, et c'est aussi
$\mathfrak{b}_{r}^{*}$. 

Supposons que l'ensemble $\mathcal{S}$ des \'{e}l\'{e}ments
$P$-stables de $\mathfrak{p}^{*}$ soit non vide. Il est clair que
c'est un ouvert de $\mathfrak{p}^{*}$. Il existe donc $f \in
\mathcal{S}$ tel que $f|\mathfrak{b} \in
\mathfrak{b}_{r}^{*}$. D'autre part, pour $\theta \in B$, on a $\theta
(f) \in \mathcal{S}$ et $\theta (f)|\mathfrak{b} = \theta
(f|\mathfrak{b})$. Il r\'{e}sulte alors de ce qui pr\'{e}c\`{e}de que
l'on peut supposer $f$ de la forme $f = \varphi_{\mathfrak{p}}(v)$,
avec $v = \lambda X_{\alpha_{2}} + u$, o\`{u} $\lambda \in \Bbbk$.

$\bullet$ Si $\lambda \ne 0$, on trouve facilement que
$\mathfrak{p}^{f} = \Bbbk x$, avec 
$$
\begin{array}{rl}
x  = & X_{-\alpha_{2}} + \mu_{1} X_{\alpha_{1}} + \mu_{2}
X_{\alpha_{3}} + \mu_{3} X_{\alpha_{4}} + \mu_{5} X_{\alpha_{1} 
 + \alpha_{2} + \alpha_{3}} \\
& + \mu_{6} X_{\alpha_{1} + \alpha_{2} +
\alpha_{4}} + \mu_{7} X_{\alpha_{2} + \alpha_{3} + \alpha_{4}} +
\mu_{8} X_{\alpha_{1} + 2 \alpha_{2} + \alpha_{3} + \alpha_{4}},
\end{array}
$$
o\`{u} $\mu_{1}, \dots , \mu_{8}$ sont des \'{e}l\'{e}ments non nuls
de $\Bbbk$ (d\'{e}pendant, par exemple, du choix d'une base de
Chevalley de $\mathfrak{g}$). 

$\bullet$ Si $\lambda = 0$, il vient $\mathfrak{p}^{f} = \Bbbk x$,
avec
$$
x = X_{-\alpha_{2}} + \mu_{1} X_{\alpha_{1} + \alpha_{2} + \alpha_{3}}
+ \mu_{2} X_{\alpha_{1} + \alpha_{2} + \alpha_{4}} + \mu_{3}
X_{\alpha_{2} + \alpha_{3} + \alpha_{4},}
$$
o\`{u} $\mu_{1}, \mu_{2}, \mu_{3}$ sont des scalaires non nuls.

Soit $h \in \mathfrak{h}$ d\'{e}fini par $\alpha_{1}(h) =
\alpha_{3}(h) = \alpha_{4}(h) = - \alpha_{2}(h) = 1$. Dans les deux
cas pr\'{e}c\'{e}dents, il vient $[h,x] = x$. On a alors
$[\mathfrak{p}, \mathfrak{p}^{f}] \cap \mathfrak{p}^{f} \ne \{ 0\}$,
ce qui contredit 1.7. Ainsi, $\mathfrak{p}^{*}$ ne contient aucun
\'{e}l\'{e}ment $P$-stable.

\end{document}